\documentclass[11pt,epsfig,epsf,amsfonts,amsmath,enumerate,amssymb,amsthm,graphicx]{article}
\usepackage{amsthm}
\usepackage{amssymb}
\usepackage{epic}
\usepackage{eepic}
\usepackage{graphics}
\usepackage{epsfig}
\usepackage{psfrag}

\topmargin by -\headsep \textheight 8.9in \oddsidemargin 0pt
\evensidemargin \oddsidemargin \marginparwidth 0.5in \textwidth
6.5in

\def\proof{\par\vskip 3pt\noindent\hbox{\bf Proof} :\quad}
\newtheorem*{thm*}{Theorem}

\newtheorem*{lem*}{Lemma}

\newtheorem{cor}{Corollary} [section]

\def\qed{\hfill\raisebox{3pt}{\fbox{\rule{0mm}{1mm}\hspace*{1mm}\rule{0mm}{1mm}}\,}
\vspace{8pt}}
\def\cadre{$$\vcenter\bgroup\advance\hsize by -8em\noindent
             \refstepcounter{equation}\ignorespaces}
\makeatletter
\def\endcadre{\egroup\eqno(\theequation)$$\global\@ignoretrue}
\makeatother

\def\Qset{\hbox{\hbox{Q\hskip-0.525em\lower-0.097ex
\hbox{\vrule height1.47ex width 0.07em}}\hskip0.50em}}
\def\Rset{\hbox{I\hskip-0.23em R}}

\newcounter{petit}

  \def\Pscr{{\cal P}} 
\def\Cscr{{\cal C}}   
   
  \def\Sscr{{\cal S}} 
   
\def\Gscr{{\cal G}}  
\def\Qset{\hbox{\hbox{Q\hskip-0.525em\lower-0.097ex
\hbox{\vrule height1.47ex width 0.07em}}\hskip0.50em}}
\def\Rset{\hbox{I\hskip-0.23em R}}

\begin{document}
		\vskip -0.93cm
\title{Color-critical Graphs and Hereditary Hypergraphs}
\author{Andr\'as Seb\H{o}  \\CNRS, Laboratoire G-SCOP,
	Univ. Grenoble Alpes}
\maketitle
\begin{abstract} A quick proof of Gallai's celebrated theorem on color-critical graphs is given from Gallai's simple,  ingenious lemma on factor-critical graphs, in terms of partitioning the vertex-set into a minimum number of  hyperedges of a hereditary hypergraph, generalizing the chromatic number.  We then show examples of applying the results to new problems and indicate the way to  algorithms and refined complexity results for all these examples at the same time.   
\end{abstract}

\section{Introduction }

Graphs and digraphs are without loops or parallel edges. Given a hypergraph $H=(V,E)$,  $(E\subseteq\Pscr(V),$ where $\Pscr(V)$ is the power-set of $V)$  we will call the elements of $V$ {\em vertices}, and those of $E$ {\em hyperedges}, $n:=|V|$.  A {\em cover} is a family $\Cscr\subseteq E$ such that $\displaystyle \cup \Cscr :=  \cup_{C\in \Cscr}C=V$. We suppose that $E$ is a cover. The minimum number of hyperedges in  a cover is denoted $\rho:=\rho(H)$.  The  {\em  hereditary closure} of $H$ is $H^h=(V,E^h)$ where $E^h:=\{X\subseteq e: e\in E\}$, and $H$ is {\em hereditary}, if $H^h=H$. 

 In this paper we study {\em hereditary} hypergraphs (sometimes also called  independence systems, or a simplicial complexes in the literature). Hyperedges of cardinality $1$ will be called {\em singletons}, and hyperedges of cardinality $2$ are called {\em edges}.   Deleting a vertex $v$ of the hypergraph $H$ results in the hypergraph $H-v=(V\setminus \{v\}, \{e\in E, v\notin E\} )$. For hereditary hypergraphs this is the same as deleting $v$ from all hyperedges. Like for coloring, for hereditary hypergraphs a minimum cover can be supposed to be a {\em partition} of $V$,  and {\em we will  suppose this}!  Indeed, a vertex contained in several hyperedges can be deleted from one of these hyperedges.  This assumption is of primary importance, since the edges and the singletons play a major role in such  partitions.       
 


 If $H$ is a hereditary hypergraph, we have $\rho(H)-1\le \rho(H-v)\le \rho(H)$; if the first inequality is satisfied with equality for all $v\in V$, we say that $H$ is {\em  critical}.

Given a    hypergraph  $H=(V,E)$, denote by $E_2$ the set of  its edges, $H_2:=(V,E_2)$.  The {\em components} of $H$ are defined as the those of $(H^h)_2$. These form a partition of $V$, and correspond to the usual hypergraph components: $H$ is {\em connected} if $(H^h)_2$ is connected. Abusing terminology, the vertex-set of a component will also be called component. 
The maximum size of a matching in a graph $G$ is denoted by $\nu(G)$. 
We prove Gallai's ingenious, simple lemma \cite{gallai:factorCritical} for self-containedness: 
\begin{lem*}
If $G=(V,E)$ is a connected graph, and $\nu(G-v)=\nu (G)$ for all $v\in V$, then $\nu(G)=\frac{n-1}2$. 
\end{lem*}
\proof Suppose for a contradiction that $M$ is a maximum matching and  $u\ne v\in V$ are not covered by $M$. Let the distance of $u$ and $v$ be minimum among all maximum matchings and two of their uncovered vertices.  Let $P\subseteq E(G)$ be a shortest path between $u$ and $v$. Clearly, $|P|\ge 2$,  otherwise the edge $uv$ could be added to $M$, contradicting the maximality of $M$. 

Let $u'$ be the neighbor of $u$ on $P$, and $M'$ a maximum matching of $G-u'$. Each connected component of the  symmetric difference $D$ of $M$ and $M'$ is the disjoint union of even paths and even circuits alternating between $M$ and $M'$. If $u$ and $u'$ are not in the same component of $D$, then interchanging the edges of $M$ in the component of $u'$, neither $u$ nor $u'$ is covered by a matching edge, leading to the same contradiction as before. 

On the other hand, if they are in the same component, the same interchange of the edges leads to a maximum matching that leaves  $u'$ and $v$ uncovered, contradicting the minimum choice of the distance between $u$ and $v$. 
\qed 

Gallai's theorem on color-critical graphs \cite{gallai:colorCritical} is a beautiful statement but its original proof was rather complicated, essentially more difficult than the above lemma on factor-critical graphs \cite{gallai:factorCritical}. Stehl\'\i k \cite{STCC} gave a simpler proof. We show here  that the  generalization to hereditary hypergraphs can be shortly reduced to Gallai's Lemma (Section~\ref{sec:main}), in addition giving rise  to a wide range of known and new examples (Section~\ref{sec:a})\footnote{The Theorem below and its proof have been included in a  more complex framework of an unpublished manuscript \cite{AM}. Several occurrences of old and recent, direct special cases of hereditary hypergraphs make it useful to provide an exclusive, short presentation of this general theorem, with some examples of hereditary hypergraphs of interest.}, to algorithms, and clarifications of their complexity issues. 

\section{Theorem and Proof}\label{sec:main}



\begin{thm*}
	In a connected, hereditary, critical hypergraph $\rho\le \frac{n+1}{2}$. Furthermore, either the inequality is strict and there is a minimum cover without singleton, or the equality holds, and there are minimum covers with only edges and exactly one singleton that can be any vertex.   
	
\end{thm*}

\proof  For $n\le 1$ the statement is obvious, so  suppose $H$ is critical,  $n\ge 2$. Then for all $v\in V$ there exists a minimum cover containing $\{v\}$: indeed,  adding  $\{v\}$ to a minimum cover of $H-v$,  we get a minimum cover of $H$. Consider a minimum cover of $H-v$ (partitioning $V\setminus \{v\}$),   {\em maximizing $C_v:=\cup\,\Cscr_v$, where $\Cscr_v$ is the set of its non-singleton elements}. Clearly, $\rho= |\Cscr_v|+ |V\setminus C_v|$.
	
\medskip\noindent
{\bf Claim~1}.   For all $u, v\in V$ and each component $C\subseteq V$ of $\Cscr_u\cup\Cscr_v:\, |C\cap C_u|=|C\cap C_v|.$ 

\smallskip
Indeed, if $k_u$ and $k_v$ are the number of hyperedges in this component of $\Cscr_u$ and $\Cscr_v$ respectively, then $k_u  + |C\setminus C_u|= k_v  + |C\setminus C_v|$ for  if say $k_u  + |C\setminus C_u| > k_v  + |C\setminus C_v|$, then

in the minimum cover consisting of $\Cscr_u$ and $|V\setminus C_u|$  singletons, 

replace the hyperedges in $C$, that is,    $k_u$ hyperedges of $\Cscr_u$   and $|C\setminus C_u|$ singletons, 

by the $k_v$ hyperedges of $\Cscr_v$ in $C$, and $|C\setminus C_v|$ singletons, leading to a cover of size 
\[\rho+ k_v - k_u  + |C\setminus C_v| -  |C\setminus C_u| < \rho,\]
a contradiction. But then the proven equality implies that the same replacement -- in either directions --  of the hyperedges   lead to a minimum cover, increasing the size of $C_u$ if say $|C\cap C_u|<|C\cap C_v|$, and this contradiction with the choice of $C_u$ proves the claim. 

\medskip\noindent
{\bf Claim~2}. If each minimum cover of $H$ contains a singleton, then  $\Cscr_v\subseteq E_2$ for all $v\in V$. 

\smallskip

\smallskip
Let $u\in C_v$ be arbitrary, and let us prove that $|e_u|=2$ for the hyperedge $e_u$,  $u\in e_u\in\Cscr_v$. Since  $u\in C_v\setminus C_u$, by Claim~1,  the component $C$ of $\Cscr_u\cup\Cscr_v$ containing $u$ also contains a vertex $v_0\in C_u\setminus C_v$. Let $P$ be a shortest path between $v_0$ and $u$ in  $\Cscr_u\cup\Cscr_v$ (in the connected component $C$). Let $v_0, v_1,v_2\ldots$ be the vertices of $P\subseteq E$, in fact $P\subseteq E_2$, in this order, necessarily alternating between subsets of hyperedges in $\Cscr_u$ and subsets of hyperedges in  $\Cscr_v$. We prove by induction on $|P|$  {\em the assertion that  the latter subsets (of hyperedges in $\Cscr_v$) are in fact  in $\Cscr_v$}: 

Note first that  $\{v_1,v_2\}\in\Cscr_v$, because if it was a subset of an edge $e\in\Cscr_v$, $|e|\ge 3$, then replacing $\{v_0\}$ and $e$ in the minimum cover $\Cscr_v\cup\{\{v'\}: v'\in V\setminus C_v\}$ by $\{v_0,v_1\}$, $e\setminus \{v_1\}$, we get a minimum cover, where the hyperedges of size at least two cover $C_v\cup {v_0}$ contradicting the definition of $C_v$, provided $v\ne v_0$. If $v=v_0$, $C_{v'}:=C_v\cup {v_0}$ can  occur in Claim~2 choosing any $v'\in V\setminus (C_v\cup {v_0})$; $v'$ exists, since $V\setminus (C_v\cup {v_0})\ne\emptyset$ because of the condition of Claim~2. 

This proves the assertion for $|P|=2$. 
Let  $\Cscr_{v_2}:=(\Cscr_v \setminus  \{v_1,v_2\}) \cup \{v_0,v_1\}$, and $P':=P-\{v_1,v_2\}$; $\Cscr_{v_2}$ is a minimum cover of $H-v_2$ maximizing the union of non-singletons,  and $|P'|<|P|$. Now the induction hypothesis finishes the proof of the assertion and of  Claim~2.  

\medskip
To finish the proof note first that a minimum cover without singleton implies $\rho\le\frac n2$ and we are done. 
Otherwise, Claim~2 can be applied,  and $\rho=\frac{|C_v|}2+ |V\setminus C_v|$  follows for all $v\in V$.  This formula also shows that a larger matching $\Cscr'_v$ would provide a smaller cover. So $\Cscr_v$ is a maximum matching of $H_2$ and does not cover $v$, so $\nu(H_2-v)=\nu(H_2)$ for all $v\in V$. The connectivity of $H$ means by definition that $H_2$ is connected, so the conditions of Gallai's Lemma are satisfied for $H_2$: $H_2$ is factor-critical, and $\{v\}$  $(v\in V)$ with  a perfect matching of $H_2-v$ provide a cover of size  $1+\frac{n-1}{2}=\frac{n+1}{2}$. 
\qed

Let us restate the inequality of the Theorem so that it directly contains the formulation of \cite{gallai:colorCritical}: 

\begin{cor}\label{cor:Gallai}
A  hereditary hypergraph with $n\le 2(\rho -1)$ is either not critical, or not connected.
\end{cor}











\section{Examples,  Algorithms and Conclusions}\label{sec:a}
 
In this section we show some examples applying the results to particular hypergraphs. Any hereditary hypergraph is an example, so we cannot seek for completeness,  but we try to show how the specialization works. An important surprise   is that it turned out that the role of larger hyperedges is secondary, {\em $H_2$ plays the main role:  the covers appearing in the Theorem consist only of edges; in the corollaries the components and connectivity depend only on $H_2$.} 

\begin{cor}\label{cor:concrete}
	Let $H=(V,E)$ be a  hereditary hypergraph with $|V|\le 2(\rho(H) -1)$. Then
	either there exists $v\in V$ so that $\rho(H-v)=\rho(H)-1$, 
	 or $H$ is not connected, that is, there exists a partition $\{V_1, V_2\}$ of $V$ so that $\{v_1,v_2\}\notin E$ for all $v_1\in V_1$, $v_2\in V_2$. \qed
\end{cor}


\noindent
{\bf 3.1 Hereditary hypergraphs from graphs}
\smallskip

Let $G=(V(G),E(G))$ be either an undirected graph or a digraph, the context always determines the current meaning, and we  define hereditary hypergraphs on $V(G)$. The more deeply the hypergraphs are related to $G$, the more interesting the results.
 Fix a (not necessarily finite) set $\Gscr$ of (di)graphs  and for each (di)graph $G$, let $H(G,\Gscr):=(V(G), F)$, where
  \[F:=\{U\subseteq V(G): \hbox{$U$ induces in $G$ a graph without any induced subgraph in $\Gscr$}\}.\]
When are the Theorem or its corollaries meaningful or even interesting for $H(G,\Gscr)$? 
 
  If neither of the $2$-vertex graphs are in $\Gscr$, hypergraphs  $H(G,\Gscr)$ are connected for every graph $G$,  and our Theorem and its corollaries are trivial. On $2$ vertices there are two undirected graphs:  one without, and one  with an edge between the two vertices. If the only graph of $\Gscr$  on two vertices is the edge-less graph,  $H(G,\Gscr)$ consists of cliques of $G$;  
  if it is the edge on two vertices, $H(G,\Gscr)$ consists of stable sets of $G$.
  In turn, according to Corollary~\ref{cor:concrete}, in the former case the disconnectivity of $H(G,\Gscr)$ means the disconnectivity of  $G$,  and in the latter case it means the disconnectivity of the complement of $G$. In these cases, the Theorem specializes to  Gallai's theorem.
  
  It is easy to see that in these cases the only  possibility to add to $\Gscr$  more graphs is to add a clique (or stable set) of given size. Then for some $k\in\mathbb{N}$, $H(G,\Gscr)$ is the family  of cliques or stable sets of size at most $k$ on $V$. For $k\ge 2$ the Theorem applies without change and it is then about coloring with at most $k$ vertices of each color. 
  
  \medskip\noindent
  	{\bf 3.2 Hereditary hypergraphs from digraphs}
  
  \smallskip
    Similarly, for  digraphs one of the subgraphs on two vertices has to be excluded: there are now three digraphs on two vertices: with or without an arc  as in the undirected case, or with an arc in both directions ($2$-cycle). If $\Gscr$ contains only the latter, we also do not get anything new: keeping only  arcs in both directions as an undirected edge  we reduce the problem to Gallai's colorings in undirected graphs. However, if there are some other graphs in $\Gscr$ we have three interesting special cases: cliques, stable sets (Gallai),  a third case we discuss below as also cases from multigraphs. 
\begin{cor}\label{cor:graphs]}
    	Let $G$ be a graph, $\Gscr$ a set of graphs,  $H:=H(G,\Gscr)$, and $|V(G)|\le 2(\rho(H) -1)$. Then
    	either there exists $v\in V$ so that $\rho(H-v)=\rho(H)-1$, 
    	or $H$ is not connected, that is, there exists a partition $\{V_1, V_2\}$ of $V$ so that $\{v_1,v_2\}$ for all $v_1\in V_1$, $v_2\in V_2$ induces a graph in $\Gscr$.\qed
    \end{cor}
      As argued before the corollary, the interesting cases are when {\em the unique graph on two vertices of $\Gscr$ is an edge, a non-edge or a $2$-cycle}, and in the last case there are many possibilities to exclude further induced subgraphs. For instance we can include in $\Gscr$ $3$-cycles and all graphs on $4$ vertices having $4$-cycles. Actually  an arbitrary subset of graphs having directed cycles, or the set of all such graphs can be contained in $\Gscr$, and will not make any change in the relevant critical graphs (as compared to including only $3$- and $4$-cycles, no larger hyeredge plays a role).  Corollary~\ref{cor:graphs]} holds, and partitioning into hyperedges of $H(G,\Gscr)$ means then partitioning into vertex-sets that induce acyclic digraphs: this is ``digraph coloring'', 
     for which Corollary~\ref{cor:graphs]} was asked in \cite{BS}.  (The Theorem has then already been proved, see footnote~1. Stehl\'\i k \cite{M} missed its specialization to acyclic induced subgraphs, and answered \cite{BS} using the Edmonds-Gallai structure theorem.)
 
  
  \medskip\noindent
  {\bf 3.3 More Examples}
  
  \smallskip

Clearly, common hyperedges of an arbitrary number of hereditary hypergraphs on the same set of vertices form a hereditary hypergraph. If all of them arise as stable sets of graphs, the intersection will be just the stable-set-hypergraph of the graph which is the union of the considered graphs. However, if the considered hypergraphs arise  in different ways, the intersections may provide nontrivial new cases, if the role of the edges is kept in mind.

 More generally, a {\em stable-set} in a  (not necessarily hereditary) hypergraph $H=(V,E)$ is a set $S\subseteq V$ so that $S$ does not contain any $e\in E$. (Independent sets of matroids are those that do not contain a hyperedge from the circuit-hypergraph.) The family $\Sscr$ of all stable sets is obviously a hereditary family; $S\in \Sscr$, if and only if $V\setminus S$ is a {\em transversal} or blocker of the hyperedges; the family  of transversals is an upper hereditary hypergraph, another source of examples: 
 
 In {\em upper hereditary} hypergraphs the supersets of hyperedges are also hyperedges.
 The {\em dual}  of $H=(V,E)$ is $H^d:=(V,E^d)$, where $E^d:=\{V\setminus e: e\in E\}$. The dual of a hereditary hypergraph is upper hereditary and vice versa, generating more examples; $(H^d)^d=H$. Each example of  upper hereditary hypergraphs provides an example of hereditary hypergraphs, and vice versa.  Upper hereditary hypergraphs arise for instance from vertex-sets of graphs that do contain one of a fixed set of graphs as induced subgraphs; being non-planar or non-bipartite is a special case.
 
 In multi (di)graphs $G$ with for instance edge-multipicities $z:E(G)\rightarrow \Rset$ and $\lambda\in\Rset$ we may consider the hereditary hypergraph $\{U\subseteq V(G): \hbox{sum of $z(e)$ on the edges induced by $U$}\le \lambda\}$, when  Corollary~\ref{cor:graphs]} is again meaningful. The upper bound can be replaced by any monoton function of $z(e)$ and the graph, combined with vertex multiplicities or edge- and vertex-colored graphs, $\ldots$ 
 
 \noindent
{\bf 3.4  Algorithms and Complexity}
\smallskip

 
The focus of the examples of the previous subsections was the Theorem. Algorithmic and complexity questions are less ``choosy'' and become meaningful and nontrivial for more examples.


Once in a while questions about  particular, critical hereditary hypergraphs are raised anew, sometimes as open problems like in \cite{BS} about partitioning the vertex-set into acyclic digraphs. How can the NP-hard covering participate in well-characterizing minmax theorems? The discussion of this question is beyond the possibilities of this note. 
This will be layed out in forthcoming papers, we mention the key to the solution only shortly:

It is NP-hard to compute $\rho_H$ and in hereditary hypergraphs $H$ it is not easier, since taking the hereditary closure does not affect $\rho$! 
The covering problem for the hereditary closure of $3$-uniform hypergraphs contains the $3$-dimensional matching problem \cite{GJ}, and is therefore NP-hard even if the hyperedges of the hereditary hypergraph are given explicitly, and their number is polynomial in the input size. Indeed,  $\rho = n/3$  if and only if there exists a partition into triangles. 

However, the maximum $\mu$ of the number of vertices covered by non-singletons in a cover of a hereditary hypergraph can be maximized in polynomial time, and the vertex-weighted generalization can also be solved!  It can be seen that this maximum does not change if we write here ``minimum cover'' instead of ``cover''. This allows to handle with well characterizing minmax theorems and in polynomial time some aspects of minimum covers \cite{T}, for which  results  of  Bouchet \cite{B}, Cornu\'ejols, Hartvigsen and Pulleyblank \cite{CP}, \cite{CHP} play an enlightening role. 


  \medskip\noindent
{\bf 3.5 Conclusion}:

We tried to show by the Theorem and  multiple examples how results on graph colorings may be extended to covers in hypergraphs. We continue this work with minmax and structure theorems, develop algorithms at the general level of hereditary hypergraphs, and show more applications and connections between various problems \cite{T}, \cite{AM}. 
 We hope the reader will also have the reflex of using hereditary hypergraphs when a new special case is coming up! 


\small
\begin{thebibliography}{0}	
\bibitem{B} A. Bouchet, Greedy Algorithm and Symmetric Matroids, Math. Prog. 38, 147--159 (1987). 
\bibitem{BS} J.~Bang-Jensen, T.~Bellitto, T.~Schweser, and M.~Stiebitz, Haj\'os and Ore constructions
for digraphs, arXiv:1908.0.
\bibitem{CP} G.~Cornu\'ejols, W.~R.~Pulleyblank, Critical Graphs, Matchings and Tours or a Hierarchy of Relaxations for the Travelling Salesman Problem, {Combinatorica}, {\bf 3(1)} (1983), 36--52. 
\bibitem{CHP} G.~Cornu\'ejols, D.~Hartvigsen, W.~~R~.Pulleyblank, Packing Subgraphs in a Graph, {Operations Research letters}, Volume 1, Number 4 (1982)
\bibitem{gallai:factorCritical}
T.~Gallai.
\newblock Neuer {B}eweis eines {T}utte'schen {S}atzes.
\newblock {\em A Magyar Tudom{\'a}nyos Akad{\'e}mia --- Matematikai Kutat{\'o}
  Int{\'e}zet{\'e}nek K{\" o}zlem{\'e}nyei}, 8:135--139, 1963.
\bibitem{gallai:colorCritical}
T.~Gallai.
\newblock Kritische Graphen II,
\newblock {\em A Magyar Tudom{\'a}nyos Akad{\'e}mia --- Matematikai Kutat{\'o}
  Int{\'e}zet{\'e}nek K{\" o}zlem{\'e}nyei}, 8:373--395, 1963.
  \bibitem{GJ}
Garey, Johnson, Computers and Intractability, W.~H.~Freeman, San Francisco (1979)
\bibitem{T} A.~Seb\H{o}, Minmax theorems and Algorithms for Hereditary Hypergraphs, in preparation. 
\bibitem{AM} A.~Seb\H{o},  M.~ Stehl\'\i k, ``Matching and Covering in Hereditary Hypergraphs'', manuscript 2009-2017, in preparation.
 \bibitem{STCC} M.~ Stehl\'\i k, {\em Critical graphs with connected complements}, J. of Comb. Theory, Series B 89 (2003) 189--194.
 \bibitem{M}  M.~ Stehl\'\i k, Private communication (October 2019).






\end{thebibliography}
\end{document}